\documentclass[10pt]{article}
\usepackage{amssymb}
\usepackage{amsmath}
\usepackage{epsfig}

\normalsize

\newtheorem{theorem}{Theorem}[section]
\newtheorem{lemma}[theorem]{Lemma}

\newtheorem{corollary}[theorem]{Corollary}

\newcommand{\boproof}{\noindent {\bf Proof. }}
\newcommand{\eoproof}{\hspace*{\fill} $\square$ \vspace{5pt}}
\newcommand{\red}{\sqsubseteq}

\newcommand{\R}{\mathbb R}
\newcommand{\Q}{\mathbb Q}
\newcommand{\Z}{\mathbb Z}

\newcommand{\1}{{\bf 1}}

\newcommand{\IP}{{({\text{IP}})}}

\newcommand{\AIP}{{({\text{AIP}})}}
\newcommand{\CPP}{{({\text{CIP}})}}

\newcommand{\GIP}{{{\cal G}_{\text{IP}}}}

\newcommand{\GCPP}{{{\cal H}_{\text{CIP}}}}

\begin{document}
\title{Test Sets for Integer Programs with $\Z$-Convex Objective}
\author{Raymond Hemmecke, UC Davis\\raymond@hemmecke.de}
\date{\today}
\maketitle

\begin{abstract}
In this paper we extend test set based augmentation methods for
integer linear programs to programs with more general convex
objective functions. We show existence and computability of finite
test sets for these wider problem classes by providing an explicit
relationship to Graver bases. One candidate where this new
approach may turn out fruitful is the Quadratic Assignment
Problem.
\end{abstract}


\section{Introduction} \label{Introduction}
Integer linear optimization problems
\[
\IP_{c,b}:\qquad \min\{c^\intercal z:Az=b,z\in\Z_+^n\},
\]
appear in many practical applications. One way to solve such a
problem is to start with a feasible solution $z_0$ and to replace
it by another feasible solution $z_0-v$ with smaller objective
value $c^\intercal (z_0-v)$, as long as we find such a vector
$v\in\Z^n$ that improves the current feasible solution. If the
problem is solvable, that is in particular if it is bounded, this
augmentation process has to stop (with an optimal solution).

The key step in this algorithmic scheme, besides finding an
initial feasible solution, is to find improving vectors. Universal
test sets, which depend only on the problem matrix $A$, provide
such vectors for any given $c$ and $b$ and for any non-optimal
feasible solution $z_0$ of $\IP_{c,b}$. Note that universal test
sets can in fact also be used to find an initial feasible solution
$z_0$ \cite{Hemmecke:PSP}. For a survey on all currently known
test sets for $\IP_{c,b}$ see \cite{Weismantel:98}.

Graver \cite{Graver:75} was the first to introduce a finite
universal test set. The {\bf Graver basis} $\GIP(A)$, or {\bf
Graver test set}, associated to $A$ consists of all $\red$-minimal
{\it non-zero}\/ solutions to $Az=0$, where for $u,v\in\Z^n$ we
say that $u\red v$ if $u^{(j)}v^{(j)}\geq 0$ and $|u^{(j)}|\leq
|v^{(j)}|$ for all components $j=1,\ldots,n$, that is, if $u$
belongs to the same orthant as $v$ and its components are not
greater in absolute value than the corresponding components of
$v$.

{\bf Example 1.} Consider the problem
\[
\min\{x+y:x,y\in\Z_+\}.
\]
The Graver test set associated to the problem matrix $A=0$ is
$\{\pm (1,0),\pm (0,1)\}$. As one can easily check, already the
subset $\{(1,0),(0,1)\}$ provides an improving direction to any
non-optimal solution of this particular problem instance. Thus,
with the help of $\{(1,0),(0,1)\}$, we can augment any given
feasible solution to the (in this case unique) optimal solution
$(0,0)$. \eoproof

Intrinsic to the proofs that there do exist finite (universal)
test sets for $\IP_{c,b}$ and that they do indeed provide an
improving direction to any non-optimal feasible solution, is the
fact that both the objective function and the constraints are
linear. Now let us observe what happens with a non-linear
objective function.

{\bf Example 2.} Consider the problem
\[
\min\{(x+y)^2+4(x-y)^2:x,y\in\Z_+\}.
\]
As again $A=0$, the corresponding Graver basis is $\GIP(A)=\{\pm
(1,0),\pm (0,1)\}$. However, this universal test set for the
integer {\it linear} program $\IP_{c,b}$ does not provide an
improving direction to any non-optimal feasible solution for the
{\it quadratic} problem given above:

Clearly, $(0,0)$ is again the unique optimal solution with
objective value $0$. Now consider the point $(1,1)$ with objective
value $4$. There are $4$ points reachable from $(1,1)$ via the
directions given by $\GIP(A)$: $(1,0)$ and $(0,1)$, both with
objective values $5$, and $(2,1)$ and $(1,2)$, both with objective
values $13$. Therefore, in order to reach the optimum $(0,0)$ from
$(1,1)$, additional vectors are needed in the test set.

As we will see below, the set $\GIP(A)\cup\{\pm (1,1)\}$ provides
improving directions to any non-optimal solution of the above
quadratic problem. Moreover, this property remains true even if we
change the objective function in a certain way. (For details see
below.) For example, with the directions from $\GIP(A)\cup\{\pm
(1,1)\}$ we can also find the optimum of the following program:
\[
\min\{e^{|x+y-3|}+4(x-y+2)^{6}+2x-y:x,y\in\Z_+\} \vspace{-25pt}
\]
\eoproof

In this paper we relieve the restriction to linear objective
functions and employ test set methods for the solution of integer
optimization problems
\[
\CPP_{f,b}:\qquad \min\{f(z):Az=b,z\in\Z_+^n\},
\]
where $A\in\Z^{d\times n}$, $b\in\Z^d$, and where
\[
f(z):=\sum_{i=1}^sf_i(c_i^\intercal z+c_{i,0})+c^\intercal z.
\]
Herein, $c\in\R^n$, $c_1,\ldots,c_s\in\Z^n$,
$c_{1,0},\ldots,c_{s,0}\in\Z$, and $f_i:\R\rightarrow\R$,
$i=1,\ldots,s$, are $\Z$-convex functions with minimum at $0$. We
call $g:\R\rightarrow\R$ a {\bf $\Z$-convex} function with minimum
at $\alpha\in\Z$, if the function $g(x+1)-g(x)$ is increasing on
$x\in\Z$ and if $g(x+1)-g(x)\leq 0$ for all integers $x<\alpha$
and $g(x+1)-g(x)\geq 0$ for all integers $x\geq\alpha$. Clearly,
these three conditions imply that $x=\alpha$ is a minimum of $g$
over $\Z$. We will, however, restrict our attention to $\Z$-convex
functions with minimum at $0$. This is no restriction, since we
can transform any $\Z$-convex function $g$ with minimum at
$\alpha$ to one with minimum at $0$ by considering
$\bar{g}(x)=g(x+\alpha)$ instead.

The problem type $\CPP_{f,b}$ includes for example linear integer
programs for $f_1=\cdots=f_s=0$, or quadratic integer programs for
$f_i(x)=x^2$. However, one could apply our approach also to more
exotic functions as $f_i(x)=|x|$ or $f_i(x)=-x$ for $x\leq 0$ and
$f_i(x)=e^x$ for $x>0$, that is, the functions $f_i$ considered as
functions from $\R$ to $\R$ need not be continuous.

Our main result is the following.

\begin{theorem} \label{Main Theorem}
Let $A\in\Z^{d\times n}$ and $c_1,\ldots,c_s\in\Z^n$ be given.
Denote by $C$ the $s\times n$ matrix whose rows are formed by the
vectors $c_1^\intercal,\ldots,c_s^\intercal$. Moreover, let $I_s$
denote the $s\times s$ unit matrix. Then for any particular choice
\begin{itemize}
\item of $\Z$-convex functions $f_1,\ldots,f_s$ with minima at
$x=0$,
\item of $c_{1,0},\ldots,c_{s,0}\in\Z$, and
\item of $c\in\R^n$,
\end{itemize}
the set
\[
\GCPP(A,C):=\phi_n\left(\GIP\left(\begin{array}{cc}A&0\\C&I_s\\
\end{array}\right)\right)
\]
provides an improving direction to any non-optimal feasible
solution of the problem $\CPP_{f,b}$. Herein, $\phi_n$ defines the
projection of a vector onto its first $n$ components, and for a
set $G$ of vectors $\phi_n(G)$ denotes the set of images of
elements in $G$ under $\phi_n$.
\end{theorem}

Trivially, $\GIP(A)\subseteq\GCPP(A,C)$ for any matrix $C$.
However, as we have seen in Example $2$, this inclusion can be
strict.

For $f_1=\cdots=f_s=0$, we simply obtain $\GCPP(A,0)=\GIP(A)$ as a
(universal) test set for $\IP_{c,b}$. But, as the following
example shows, the set $\GIP(A)$ gives improving directions even
for a far bigger problem class.

{\bf Example 3.} Consider the family of problems $\CPP_{f,b}$
where $c_1,\ldots,c_n$ are the unit vectors in $\R^n$, that is,
\[
\min\{f(z):Az=b,z\in\Z_+^n\}
\]
with
\[
f(z):=\sum_{i=1}^sf_i(z_i+c_{i,0})+c^\intercal z.
\]
As $C=I_n$, we need to compute the Graver basis of the Lawrence
lifting
\[
\left(\begin{array}{cc}A&0\\I_n&I_n\\
\end{array}\right)
\]
of $A$. Since all elements in the kernel of this Lawrence lifting
have the form $(u,-u)$ and since $(v,-v)\red (u,-u)$ in $\Z^{2n}$
if and only if $v\red u$ in $\Z^n$, this Graver basis is simply
$\{(u,-u):u\in\GIP(A)\}$. Thus, $\GCPP(A,I_n)=\GIP(A)$, showing
that the set $\GIP(A)$ is also a test set for this bigger problem
class where $A$ is kept fixed and the remaining problem data is
allowed to vary. \eoproof

Although test set based methods are not yet proven to be
successful in practice, there is renewed hope from recent work on
generating functions \cite{Barvinok+Woods:2003,LattE:2003}, in
which it is proved that in fixed dimension any given problem
$\IP_{c,b}$ can be solved via test sets in time polynomial in the
input data. It would be an interesting research project to
generalize this complexity result to certain classes of functions
$f_i$, for example to $f_i(x)=\alpha_ix^{2\gamma_i}$ with
$\alpha_i>0$ and $\gamma_i\in\Z_+$.

The remainder of this paper is structured as follows: In Section
\ref{Section: Quadratic Programs} we show that our test set
approach can be applied to convex quadratic optimization problems,
of which the Quadratic Assignment Problem (QAP) is probably the
most famous example. Finally, in Section \ref{Section: Proof of
Main Theorem} we prove our main theorem, Theorem \ref{Main
Theorem}.

\section{Quadratic Programs}\label{Section: Quadratic Programs}
In this section we deal with the special case of convex quadratic
optimization problems
\[
\min\{z^\intercal Qz+c^\intercal z: Az=b, z\in\Z_+^n\},
\]
where $Q$ is a symmetric, positive semi-definite matrix with only
rational entries. These problems can be solved by the test set
approach introduced in Section \ref{Introduction}. The reason for
this is the following basic result from the theory of quadratic
forms \cite{Lam:73}.

\begin{lemma} \label{Q=U'DU general}
Let $Q\in\Q^{n\times n}$ be a symmetric matrix. Then there exist a
diagonal matrix $D\in\Q^{n\times n}$ and an invertible matrix
$U\in\Q^{n\times n}$ such that $Q=U^\intercal D U$. Moreover, each
diagonal element $d_{ii}$ of $D$ is representable by the quadratic
form $x^\intercal Qx$, that is, for all $d_{ii}$ there is some
$x_i\in\R^n$ such that $d_{ii}=x_i^\intercal Qx_i$.
\end{lemma}

\begin{corollary} \label{Q=U'DU}
Let $Q\in\Q^{n\times n}$ be a symmetric positive semi-definite
matrix. Then there exist a diagonal matrix $D\in\Q^{n\times n}$ with
only non-negative entries and an invertible matrix
$U\in\Q^{n\times n}$ such that $Q=U^\intercal D U$.
\end{corollary}

\boproof This is an immediate consequence of Lemma \ref{Q=U'DU
general}, since $d_{ii}=x_i^\intercal Qx_i\geq 0$ for all $i$ as
$Q$ is positive semi-definite. \eoproof

Thus, every convex quadratic objective function $z^\intercal Qz$
can be restated as $\sum_{i=1}^s\alpha_i(c_i^\intercal z)^2$ with
$\alpha_i>0$ and $c_i\in\Z^n$. Therefore, the test set approach
presented in Section \ref{Introduction} is applicable to these
problems with $f_i(x)=\alpha_i x_i^2$, $\alpha_i>0$. Moreover, we
should point out that $s\leq n$, that is, the Graver basis that
has to be computed for $\GCPP(A,C)$ involves at most $2n$
variables.

In the following, we will restrict our attention to quadratic
$0$-$1$ problems.

\begin{corollary} \label{Q=U'DU for 0-1}
Any quadratic $0$-$1$ optimization problem
\[
\min\{z^\intercal Qz+c^\intercal z: Az=b, z\in\{0,1\}^n\}
\]
with symmetric matrix $Q\in\Q^{n\times n}$ can be rephrased as an
equivalent problem
\[
\min\{z^\intercal\bar{Q}z+\bar{c}^\intercal z: Az=b,
z\in\{0,1\}^n\},
\]
where $\bar{Q}\in\Q^{n\times n}$ is a symmetric, positive definite
matrix.
\end{corollary}

\boproof As $z_i^2-z_i=0$ for $z\in\{0,1\}$, the given
optimization problem is equivalent to
\[
\min\{z^\intercal Qz+c^\intercal z+\lambda(z^\intercal I_nz-
\1^\intercal z): Az=b, z\in\{0,1\}^n\},
\]
where $\lambda\in\R$ denotes some fixed scalar. As for
sufficiently large $\lambda=\bar{\lambda}\in\Q_+$ the matrix
$Q+\bar{\lambda} I_n$ becomes positive definite, Lemma
\ref{Q=U'DU} can be applied, giving the result with
$\bar{Q}=Q+\bar{\lambda} I_n$ and $\bar{c}=c-\bar{\lambda}\1$.
\eoproof

Consequently, {\it any} $0$-$1$ quadratic optimization problem
\[
\min\{z^\intercal Qz+c^\intercal z: Az=b, z\in\{0,1\}^n\}
\]
can be written as
\[
\min\{\sum_{i=1}^s\alpha_i(c_i^\intercal z)^2+c^\intercal z:Az=b,
z\in\{0,1\}^n\}
\]
with $\alpha_i>0$, and therefore the test set approach presented
in Section \ref{Introduction} can be applied. However, choosing
different $\bar{\lambda}$ in the proof of Corollary \ref{Q=U'DU
for 0-1}, we get different equivalent formulations for the same
problem $\CPP_{f,b}$. But as the following example shows,
different problem formulations can lead to different test sets
$\GCPP(A,C)$ for the {\it same} problem. These sets, however, are
test sets for two {\it different problem families} of which the
given specific problem is a {\it common member}.

{\bf Example 4.} Consider the quadratic $0$-$1$ problem with $A=0$
and
\[
Q=\left(
\begin{array}{ccc}
0 & 1 & 1\\
1 & 0 & 2\\
1 & 2 & 0\\
\end{array}
\right).
\]
Since $A=0$, we need to compute the Graver basis of $(C|I_3)$. But
we have different choices for $C$. As $x_i^2=x_i$, $i=1,2,3$, we
have
\begin{eqnarray*}
x^\intercal Qx & = & 2x_1x_2+2x_1x_3+4x_2x_3\\
& = & (x_1+x_2+x_3)^2+(x_2+x_3)^2-x_1^2-2x_2^2-2x_3^2\\
& = & (x_1+x_2+x_3)^2+(x_2+x_3)^2-x_1-2x_2-2x_3\\
\end{eqnarray*}
and
\begin{eqnarray*}
x^\intercal Qx & = & 2x_1x_2+2x_1x_3+4x_2x_3\\
& = & (x_1-2x_2+x_3)^2+(3x_1+x_2+4x_3)^2+12(x_1-x_3)^2-
22x_1^2-5x_2^2-29x_3^2\\
& = & (x_1-2x_2+x_3)^2+(3x_1+x_2+4x_3)^2+12(x_1-x_3)^2-
22x_1-5x_2-29x_3.\\
\end{eqnarray*}
Therefore, the corresponding two matrices for the test set
computations are
\[
\left(
\begin{array}{c|c}
C' & I_2\\
\end{array}
\right)= \left(
\begin{array}{ccc|cc}
1 & 1 & 1 & 1 & 0\\
0 & 1 & 1 & 0 & 1\\
\end{array}
\right)
\]
and
\[
\left(
\begin{array}{cc}
C'' & I_3\\
\end{array}
\right)= \left(
\begin{array}{rrr|rrr}
1 & -2 &  1 & 1 & 0 & 0\\
3 &  1 &  4 & 0 & 1 & 0\\
1 &  0 & -1 & 0 & 0 & 1\\
\end{array}
\right).
\]
Using the software package 4ti2 \cite{4ti2}, we obtain
\begin{eqnarray*}
\GCPP(A',C') & = & \{(1,0,0),(0,1,0),(0, 0,
1),(1,-1,0),(0,1,-1),(1,0,-1)\}\\
\GCPP(A'',C'')& = & \{(1, 0,
0),(0, 1, 0),(0, 0, 1),(1, 1, 0),(1, -1, 0),
(0, 1, -1),\\
& & (0, 1, 1),(1, 0, 1),(1, 0, -1),(1, 1, -1),(1, 1, 1),(1, -1,
1)\}
\end{eqnarray*}
Note that $\GCPP(A',C')\subsetneq\GCPP(A'',C'')$. \eoproof

This gives us much freedom to rewrite particular $0$-$1$ problems,
possibly arriving at much smaller test sets for the same problem.
As the following example shows, the same phenomenon happens also
in the general (non-$0$-$1$) case.

{\bf Example 5.} Consider the problem with $A=0$ and
\[
Q=\left(
\begin{array}{ccc}
2 & 1 & 1\\
1 & 2 & 1\\
1 & 1 & 2\\
\end{array}
\right).
\]
Again, since $A=0$, we need to compute the Graver basis of
$(C|I_s)$ for some integer $s$, and as the following shows, we
have more than one choice for $C$:
\begin{eqnarray*}
x^\intercal Qx & = & 2x_1^2+2x_2^2+2x_3^2+2x_1x_2+2x_1x_3+2x_2x_3\\
& = & (x_1+x_2+x_3)^2+x_1^2+x_2^2+x_3^2\\
& = & (x_1+x_2)^2+(x_1+x_3)^2+(x_2+x_3)^2
\end{eqnarray*}
Corresponding to these two representations are the matrices
\[
\left(
\begin{array}{c|c}
C' & I_4\\
\end{array}
\right)= \left(
\begin{array}{ccc|cccc}
1 & 1 & 1 & 1 & 0 & 0 & 0\\
1 & 0 & 0 & 0 & 1 & 0 & 0\\
0 & 1 & 0 & 0 & 0 & 1 & 0\\
0 & 0 & 1 & 0 & 0 & 0 & 1\\
\end{array}
\right)
\]
and
\[
\left(
\begin{array}{c|c}
C'' & I_3\\
\end{array}
\right)= \left(
\begin{array}{ccc|ccc}
1 & 1 & 0 & 1 & 0 & 0\\
1 & 0 & 1 & 0 & 1 & 0\\
0 & 1 & 1 & 0 & 0 & 1\\
\end{array}
\right).
\]
Using 4ti2 again, we obtain
\begin{eqnarray*}
\GCPP(A',C') & = & \{(1,0,0),(0,1,0),(0,0,1),(1,-1,0),(0,1,-1),(1,0,-1)\}\\
\GCPP(A'',C'')& = & \{(1,0,0),(0,1,0),(0,0,1),(1,-1,0),(1,0,-1),(0,1,-1),\\
& & (1,1,-1),(1,-1,1),(1,-1,-1)\}
\end{eqnarray*}
Note that again, $\GCPP(A',C')\subsetneq\GCPP(A'',C'')$. \eoproof

The quadratic assignment problem
\cite{Burkard+Cela+Pardalos+Pitsoulis:1998} deals with assigning
$n$ facilities to $n$ locations such that a certain quadratic cost
function is minimized. It can be formulated as the following
problem involving permutation matrices :
\begin{eqnarray*}
\min\{\sum_{i=1}^n\sum_{j=1}^n\sum_{k=1}^n\sum_{l=1}^n
d_{ijkl}x_{ij}x_{kl} & + & \sum_{i=1}^n\sum_{j=1}^n c_{ij}x_{ij}:\\
& & \sum_{j=1}^n x_{ij} = 1, \;\;\; i\in\{1,\ldots,n\},\\
& & \sum_{i=1}^n x_{ij} = 1, \;\;\; j\in\{1,\ldots,n\},\\
& & x_{ij}\in\{0,1\}, \;\;\; i,j\in\{1,\ldots,n\}.
\end{eqnarray*}
The value $d_{ijkl}$ can be seen as costs for assigning facility
$i$ to location $j$ {\bf and} facility $k$ to location $l$,
whereas $c_{ij}$ models a fixed cost incurred by locating facility
$i$ to location $j$.

Even nowadays, QAP's of size $n>30$ (that is, with more than only
$900$ binary variables) are still considered to be computationally
extremely hard, if not intractable. One major problem in
branch-and-bound algorithms that try to solve these problems is
the lack of sharp lower bounds.

As we had seen after Corollary \ref{Q=U'DU}, our novel approach
presented in Section \ref{Introduction} reduces the question of
solving the QAP to finding a truncated Graver basis in at most
$2n^2$ variables, of which $n^2$ variables are bounded by $1$.

From a practical perspective, however, we can restrict our
attention to certain orthants to find an improving vector to a
given feasible $0$-$1$ solution. Moreover, we can use the upper
bound of $1$. Besides speeding up the computation, both constraint
types reduce drastically the number of test set vectors that could
provide an improving direction to the current solution, a very
important fact for practical applicability.

We think it to be an interesting future project to try our new
test set approach to instances from the QAPLIB \cite{QAPLIB}.
Although the software package 4ti2 \cite{4ti2} exploits both
orthant and upper bound constraints, it does not yet include a
special $0$-$1$ implementation in which special data structures
speed up the computation and save valuable memory.

\section{Proof of Main Theorem} \label{Section: Proof of Main Theorem}

In this section we prove the main theorem, Theorem \ref{Main
Theorem}, of this paper. First, we will collect some facts about
Graver bases that will turn out very useful in the final proof.
Lemma $3.23$ in \cite{Hemmecke:SIP2} states the following.

\begin{lemma} \label{IP G(A a -a)}
Let $B=\left(\begin{array}{rrr} A & a & -a\\ \end{array}\right)$
be an integer matrix such that the two columns $a$ and $-a$ differ
only by a sign. Then the Graver basis of $B$ can be constructed
from the Graver basis of $B'=\left(\begin{array}{rr} A & a
\\ \end{array}\right)$ in the following way:
\[
\GIP(B)=\{(u,v,w):vw\leq 0, (u,v-w)\in\GIP(B')\}\cup\{\pm
(0,1,1)\}.
\]
\end{lemma}
A simple corollary of this is
\begin{corollary}\label{IP G(A a a)}
Let $B=\left(\begin{array}{rrr} A & a & a\\ \end{array}\right)$ be
an integer matrix with two identical columns $a$. Then the Graver
basis of $B$ can be constructed from the Graver basis of
$B'=\left(\begin{array}{rr} A & a \\ \end{array}\right)$ in the
following way:
\[
\GIP(B)=\{(u,v,w):vw\geq 0, (u,v+w)\in\GIP(B')\}\cup\{\pm
(0,1,-1)\}.
\]
\end{corollary}

\boproof The claim follows immediately from the fact that
$(u,v,w)$ is $\red$-minimal in $\ker\left(\begin{array}{rrr} A & a
& a\\ \end{array}\right)$ if and only if $(u,v,-w)$ is
$\red$-minimal in $\ker\left(\begin{array}{rrr} A & a & -a\\
\end{array}\right)$. \eoproof

The following is an immediate consequence of Lemma \ref{IP G(A a
-a)} and of Corollary \ref{IP G(A a a)}.

\begin{lemma} \label{(A a ... a -a ... -a)}
Let $A\in\Z^{d\times n}$ and let $B=\left(\begin{array}{rrrrrrr} A
& a & \ldots & a & -a & \ldots & -a\\ \end{array}\right)$ be an
integer matrix with finitely many multiple columns $a$ and $-a$
which differ only in their signs. Then we have
$\phi_n(\GIP(B))=\phi_n\left(\GIP\left(\begin{array}{rr} A & a\\
\end{array}\right)\right)\cup\{0\}$.
\end{lemma}

\boproof The constructions in Lemma \ref{IP G(A a -a)} and in
Corollary \ref{IP G(A a a)} satisfy
$\phi_n\left(\GIP\left(\begin{array}{rrr} A & a & -a\\
\end{array}\right)\right) =\phi_n\left(\GIP\left(\begin{array}{rr} A &
a\\ \end{array}\right)\right)\cup\{0\}$ and
$\phi_n\left(\GIP\left(\begin{array}{rrr} A & a & a\\
\end{array}\right)\right)=\phi_n\left(\GIP\left(\begin{array}{rr} A & a\\
\end{array}\right)\right)\cup\{0\}$. Putting both constructions
together iteratively, we get
$\phi_n(\GIP(B))=\phi_n\left(\GIP\left(\begin{array}{rr} A & a\\
\end{array}\right)\right)\cup\{0\}$, as claimed. \eoproof

Thus, in order to compute $\phi_n(\GIP(B))$, it suffices to
compute $\phi_n\left(\GIP\left(\begin{array}{rr} A & a\\
\end{array}\right)\right)$. The following is an immediate
consequence to Lemma \ref{(A a ... a -a ... -a)}.

\begin{corollary} \label{Corollary to (A a ... a -a ... -a)}
Let $A\in\Z^{d\times n}$, $c_1,\ldots,c_s\in\Z^n$, and
$k\in\Z_{>0}$. Denote by $C$ the $s\times n$ matrix whose rows are
formed by the vectors $c_1^\intercal,\ldots,c_s^\intercal$, by the
bold letter $\1$ the vector in $R^k$ with all entries $1$, and by
$I_s$ the $s\times s$ unit matrix. Then
\[
\phi_n(\GIP(A_k))=\phi_n\left(\GIP\left(\begin{array}{cc}A&0\\C&I_s\\
\end{array}\right)\right)\cup\{0\},
\]
where
\[
A_k:=\left(
\begin{array}{ccccccccc}
A             &     &    &     &    &        &        &     &   \\
c_1^\intercal & -\1 & \1 &     &    &        &        &     &   \\
c_2^\intercal &     &    & -\1 & \1 &        &        &     &   \\
              &     &    &     &    & \ddots & \ddots &     &   \\
c_s^\intercal &     &    &     &    &        &        & -\1 & \1\\
\end{array}
\right).
\]
\end{corollary}

Before we come to the proof of our main theorem, let us prove two
more useful facts.

\begin{lemma} \label{Lemma: optimal solution for one function}
Let $g$ be a $\Z$-convex function with minimum at $0$. Then for
fixed $p\in\Z$ and for fixed $k\geq |p|$, an optimal solution to
\[
\begin{array}{lrl}
\min\{\sum\limits_{j=1}^k &
(g(j)-g(j-1))x_{i,j}+(g(-j)-g(-j+1))y_{i,j}: &\\
& p=\sum\limits_{j=1}^k x_{j}-\sum\limits_{j=1}^k y_{j}, \;\;\; x_{j},y_{j}\in\{0,1\}, & j=1,\ldots,k\},\\
\end{array}
\]
is given by
\[
\begin{array}{lcl}
x_1=\ldots=x_p=1, x_{p+1}=\ldots=x_k=y_1=\ldots=y_k=0, & \text{if} & p>0,\\
x_1=\ldots=x_k=y_1=\ldots=y_k=0, & \text{if} &p=0,\\
y_1=\ldots=y_{-p}=1, y_{-p+1}=\ldots=y_k=x_1=\ldots=x_k=0, & \text{if} & p<0.\\
\end{array}
\]
The optimal value in each of these three cases is $g(p)-g(0)$.
\end{lemma}

\boproof The case $p=0$ is trivial and the optimal objective value
is $0=g(0)-g(0)$.

Let us now consider the case $p>0$. Clearly, since $p>0$, some
$x_i$ must be positive. Suppose that in a minimal solution we have
$x_i=1$ and $y_j=1$ for some $i$ and some $j$. This cannot happen,
since by putting $x_i=0$ and $y_j=0$ we would arrive at a solution
with smaller objective value, as all coefficients in the objective
function are positive. Thus, in a minimal solution
$y_1=\ldots=y_k=0$.

Since $g$ is a $\Z$-convex function with minimum at $0$, the
coefficients $g(j)-g(j-1)$ in the objective function are
non-negative and form an increasing sequence as $j>0$ increases.
Thus, $x_1=\ldots=x_p=1$, $x_{p+1}=\ldots=x_k=0$ leads to a
minimal objective value. This value is
\[
\sum\limits_{j=1}^p (g(j)-g(j-1))=g(p)-g(0).
\]

For the case $p<0$ we conclude analogously that
$x_1=\ldots=x_k=0$. Moreover, since $g$ is a $\Z$-convex function
with minimum at $0$, the coefficients $g(-j)-g(-j+1)$ in the
objective function are non-negative and form an increasing
sequence as $j>0$ increases. As above, this implies that
$y_1=\ldots=y_{-p}=1$, $y_{-p+1}=\ldots=y_k=0$ leads to a minimal
objective value. This value is again
\[
\sum\limits_{j=1}^{-p} (g(-j)-g(-j+1))=g(p)-g(0)
\]
and the claim is proved. \eoproof

\begin{lemma} \label{Lemma: optimal solution for full function}
Let $f_1,\ldots,f_s$ be $\Z$-convex functions with minimum at $0$,
$A\in\Z^{d\times n}$, $b\in\Z^d$, $c\in\R^n$,
$c_1,\ldots,c_s\in\Z^n$, and $c_{1,0},\ldots,c_{s,0}\in\Z$ be
given. Then for fixed $z\in\Z^n$ and for fixed
$k\geq\max\{|c_i^\intercal z+c_{i,0}|,i=1,\ldots,s\}$, the optimal
value of
\[
\begin{array}{lll}
\min\{\sum\limits_{i=1}^s\sum\limits_{j=1}^k &
(f_i(j)-f_i(j-1))x_{i,j}+(f_i(-j)-f_i(-j+1))y_{i,j}+c^\intercal z: &\\
& Az=b, & z\in\Z_+^n,\\
& c_i^\intercal z+c_{i,0}=\sum\limits_{j=1}^k x_{i,j}-\sum\limits_{j=1}^k y_{i,j}, & i=1,\ldots,s,\\
& x_{i,j},y_{i,j}\in\{0,1\}, & i=1,\ldots,s,\\
& & j=1,\ldots,k\}.\\
\end{array}
\]
is $f(z)-\sum\limits_{i=1}^sf_i(0)$, where
\[
f(z):=\sum\limits_{i=1}^sf_i(c_i^\intercal z+c_{i,0})+c^\intercal
z.
\]
\end{lemma}

\boproof Since $z$ is fixed, the problem decomposes into $s$
smaller problems for which we can apply Lemma \ref{Lemma: optimal
solution for one function}. Thus, the optimal value of the given
problem is
\[
\sum_{i=1}^s[f_i(c_i^\intercal z+c_{i,0})-f_i(0)]+c^\intercal z =
\sum_{i=1}^sf_i(c_i^\intercal z+c_{i,0})+c^\intercal
z-\sum_{i=1}^sf_i(0) = f(z)-\sum_{i=1}^sf_i(0).
\]
\eoproof

Now let us finally prove our main theorem, Theorem \ref{Main
Theorem}, introduced in Section \ref{Introduction}.

\boproof In order to prove this claim, assume that we are given
$\Z$-convex functions $f_1,\ldots,f_s$ with minimum at $0$,
$c_{1,0},\ldots,c_{s,0}\in\Z$, $b\in\Z^d$, and $c\in\R^n$.
Moreover, assume that we are given a non-optimal feasible solution
$z_0$ to $Az=b$, $z\in\Z_+^n$.

The theorem is proved if we can find some vector $t\in\GCPP(A,C)$
such that $z_0-t$ is feasible and such that $f(z_0-t)<f(z_0)$. In
the following, we construct such a vector $t$.

Since we assume $z_0$ to be non-minimal, there exists some better
feasible solution $z_1$, say. Let
\[
k:=\max\{|c_i^\intercal z_0+c_{i,0}|,|c_i^\intercal
z_1+c_{i,0}|,i=1,\ldots,s\}
\]
and consider the auxiliary integer linear program
\[
\begin{array}{lll}
\AIP:\;\; \min\{\sum\limits_{i=1}^s\sum\limits_{j=1}^k &
(f_i(j)-f_i(j-1))x_{i,j}+(f_i(-j)-f_i(-j+1))y_{i,j}+c^\intercal z: &\\
& Az=b, & z\in\Z_+^n,\\
& c_i^\intercal z+c_{i,0}=\sum\limits_{j=1}^k x_{i,j}-\sum\limits_{j=1}^k y_{i,j}, & i=1,\ldots,s,\\
& x_{i,j},y_{i,j}\in\{0,1\}, & i=1,\ldots,s,\\
& & j=1,\ldots,k\}.\\
\end{array}
\]
By Lemmas \ref{Lemma: optimal solution for one function} and
\ref{Lemma: optimal solution for full function}, the minimal
values of $\AIP$ for fixed $z=z_0$ and $z=z_1$ are $f(z_0)-f_0$
and $f(z_1)-f_0$, where $f_0=\sum\limits_{i=1}^sf_i(0)$. By
$(z_0,x_0,y_0)$ and $(z_1,x_1,y_1)$ denote feasible solutions of
$\AIP$ that achieve these values.

As $f(z_0)>f(z_1)$ by assumption, we have $f(z_0)-f_0>f(z_1)-f_0$.
Thus, $(z_0,x_0,y_0)$ is a feasible solution of $\AIP$ that is not
optimal. Therefore, there must exist some vector $(t,u,v)$ in the
Graver basis associated with the problem matrix of $\AIP$ that
improves $(z_0,x_0,y_0)$. We will now show that $t\in\GCPP(A,C)$,
that $z_0-t$ is feasible for $\CPP_{f,b}$, and that
$f(z_0-t)<f(z_0)$. The claim then follows immediately.

The problem matrix associated to $\AIP$ is
\[
A_k:=\left(
\begin{array}{ccccccccc}
A             &     &    &     &    &        &        &     &   \\
c_1^\intercal & -\1 & \1 &     &    &        &        &     &   \\
c_2^\intercal &     &    & -\1 & \1 &        &        &     &   \\
              &     &    &     &    & \ddots & \ddots &     &   \\
c_s^\intercal &     &    &     &    &        &        & -\1 & \1\\
\end{array}
\right),
\]
where
\[
\phi_n(\GIP(A_k))=\phi_n\left(\GIP\left(\begin{array}{cc}A&0\\C&I_s\\
\end{array}\right)\right)\cup\{0\}=\GCPP(A,C)\cup\{0\},
\]
by Corollary \ref{Corollary to (A a ... a -a ... -a)}. Therefore,
$(t,u,v)\in\GIP(A_k)$ satisfies $t\in\GCPP(A,C)\cup\{0\}$.
Moreover, as $(z_0,x_0,y_0)-(t,u,v)$ is feasible for $\AIP$, we
must have $A(z_0-t)=b$ and $z_0-t\geq 0$, implying that $z_0-t$ is
feasible for $\CPP_{f,b}$.

It remains to show $f(z_0-t)<f(z_0)$, since this also implies
$t\neq 0$ and hence $t\in\GCPP(A,C)$.

Let $(z_0-t,x_2,y_2)$ be a feasible solution of $\AIP$ that
achieves the minimal value $f(z_0-t)-f_0$ of $\AIP$ for fixed
$z=z_0-t$, see Lemmas \ref{Lemma: optimal solution for one
function} and \ref{Lemma: optimal solution for full function} for
its existence and construction. Clearly, this minimal objective
value for fixed $z=z_0-t$ is less than or equal to the objective
value of $(z_0,x_0,y_0)-(t,u,v)$, which in turn is strictly less
than $f(z_0)-f_0$, the objective value of $(z_0,x_0,y_0)$.

Therefore, $f(z_0-t)-f_0<f(z_0)-f_0$ and consequently
$f(z_0-t)-f(z_0)$. \eoproof

{\bf Acknowledgment.} The author would like to thank Kristen Nairn
for many helpful comments on this paper.

\end{document}